# ITERATED BROWNIAN MOTION IN AN OPEN SET

By R. Dante DeBlassie

*Texas A&M University*

Suppose a solid has a crack filled with a gas. If the crack reaches the surrounding medium, how long does it take the gas to diffuse out of the crack? Iterated Brownian motion serves as a model for diffusion in a crack. If $\tau$ is the first exit time of iterated Brownian motion from the solid, then $P(\tau > t)$ can be viewed as a measurement of the amount of contaminant left in the crack at time $t$. We determine the large time asymptotics of $P(\tau > t)$ for both bounded and unbounded sets. We also discuss a strange connection between iterated Brownian motion and the parabolic operator $\frac{1}{8}\Delta^2 - \frac{\partial}{\partial t}$.

**1. Introduction.** Suppose an infinite slab with finite thickness is cracked. The crack is filled with a gas and the crack reaches the surrounding medium. How long does it take the gas to diffuse out of the crack? Burdzy and Khoshnevisan [14] show a reasonable model of diffusion in a crack is given by iterated Brownian motion. In essence, the crack is modeled as the path of two-sided Brownian motion. In analogy with ordinary Brownian motion and diffusion, if $\tau$ is the first exit time of iterated Brownian motion from the solid, then $P(\tau > t)$ provides a measure of the amount of contaminant left after time $t$. This is the object of our study.

To define iterated Brownian motion $Z_t$ started at $z \in \mathbb{R}$, let $X_t^+, X_t^-$ and $Y_t$ be independent one-dimensional Brownian motions, all started at 0. Two-sided Brownian motion is defined to be

$$X_t = \begin{cases} X_t^+, & t \geq 0, \\ X_{(-t)}^-, & t < 0. \end{cases}$$

Then iterated Brownian motion started at $z \in \mathbb{R}$ is

$$Z_t = z + X(Y_t), \qquad t \geq 0.$$

In $\mathbb{R}^n$ one requires $X^\pm$ to be independent $n$-dimensional Brownian motions. Observe there is a choice here: ours or the requirement that the components of $n$-dimensional iterated Brownian motion be independent one-dimensional iterated Brownian motions. In fact, the latter process is the subject of the article [48] by Xiao. Our choice is motivated by a connection with $\frac{1}{8}\Delta^2 - \frac{\partial}{\partial t}$ and the interpretation of Burdzy and Khoshnevisan [14] as diffusion in a Brownian crack.

Recently this process has been the subject of many articles. Limit theorems are studied in [4, 17, 21, 29, 32, 46]; path properties and local time in [8, 11, 18, 30, 31, 48]; level sets in [13] and [28]. Burdzy [12] studied the variation of iterated Brownian motion. These results were extended by Khoshnevisan and Lewis [33] in their development of a stochastic integral for iterated Brownian motion.

Although iterated Brownian motion is not a Markov process (the Chapman–Kolmogorov equation is not valid), a variant of it is known to have a connection with the parabolic operator $\frac{1}{8}\Delta^2 - \frac{\partial}{\partial t}$ [23].

The purpose of this article is two-fold: First, we establish a strange connection between iterated Brownian motion and a certain partial differential equation involving $\frac{1}{8}\Delta^2 - \frac{\partial}{\partial t}$. If this connection is analogous to the Brownian case, there might be a way to study the lifetime of iterated Brownian motion in an open set by solving an initial-boundary value problem. We show this is not true in general. The second purpose of the article is to determine asymptotics of the distribution of the lifetime using another method. Now we write detailed statements.

For typographical simplicity, write

$$p_n(t,u) = \frac{1}{(2\pi t)^{n/2}} \exp\left(-\frac{u^2}{2t}\right).$$

If $f \in C_b(\mathbb{R}^n)$, the space of bounded and continuous real-valued functions on $\mathbb{R}^n$, set

$$u(t,x) = E_x(f(Z_t)),$$

where $E_x$ denotes expectation associated with $Z_0 = x$. Then

$$
\begin{aligned}
u(t,x) &= E\Big(f(x + X^+(0 \vee Y_t) + X^-(0 \vee (-Y_t)))\Big) \\
&= \int_{-\infty}^{\infty} E\Big(f(x + X^+(0 \vee y) + X^-(0 \vee (-y)))\Big) p_1(t,y)\,dy \\
&= 2\int_0^{\infty} \int_{\mathbb{R}^n} f(w) p_n(y, |w-x|) p_1(t,y)\,dw\,dy.
\end{aligned}
$$
(1.1)

If we proceed *formally* and differentiate under the integral with respect to $x$, using that $\Delta p_n(t,x) = 2\frac{\partial}{\partial t}p_n(t,x)$, we get

$$
\begin{aligned}
\Delta^2 u &= 2\int_0^{\infty} \int_{\mathbb{R}^n} f(w) \Delta_x \Delta_x(p_n(y, |w-x|)) p_1(t,y)\,dw\,dy \\
&= 8\int_0^{\infty} \int_{\mathbb{R}^n} f(w) \left[\frac{\partial^2}{\partial y^2} p_n(y, |w-x|)\right] p_1(t,y)\,dw\,dy.
\end{aligned}
$$

Next, change the order of integration, integrate by parts twice and then reverse the order of integration again to get

$$\Delta^2 u = 8 \int_0^\infty \int_{\mathbb{R}^n} f(w) p_n(y, |w-x|) \frac{\partial^2}{\partial y^2} p_1(t, y) \, dw \, dy$$
$$= 16 \frac{\partial u}{\partial t}.$$

Thus, there is hope a slight variation of Funaki's result carries through for our version of iterated Brownian motion. This is false. Differentiation under the integral is not possible and this leads to the following strange connection.

THEOREM 1.1. *Let $f \in C_b^\infty(\mathbb{R}^n)$. Then the function $u(t, x) = E_x[f(Z_t)]$ solves the Cauchy problem*

$$u \in C^\infty((0, \infty) \times \mathbb{R}^n) \cap C_b([0, \infty) \times \mathbb{R}^n),$$

(1.2) $\quad \left(\frac{1}{8}\Delta^2 - \frac{\partial}{\partial t}\right) u(t, x) = -\frac{1}{2} \frac{1}{\sqrt{2\pi t}} \Delta f(x) \qquad on\ (0, \infty) \times \mathbb{R}^n,$

$$u(0, x) = f(x).$$

After submitting this article, an elegant variant of Theorem 1.1 due to Allouba and Zheng [2] appeared in print. See Theorem 0.1 as well as other interesting results in that paper. The theorem generalizes our Theorem 1.1 in two ways: First, the Laplacian $\Delta$ can be replaced by the generator $\mathcal{A}$ of a continuous Markov process. Thus, iterated Brownian motion can be replaced by a more general iterated Markov process. Second, the function $f$ need only be a bounded element of the domain of $\mathcal{A}$ with bounded Hölder continuous second-order partial derivatives.

Let $B_t$ be $n$-dimensional Brownian motion and $D \subseteq \mathbb{R}^n$ a reasonable open set. For $\tau_D(B) = \inf\{t \geq 0 : B_t \notin D\}$, it is well known that the function

$$v(t, x) = E_x[g(B_t) I(\tau_D(B) > t)]$$

solves the initial-boundary value problem

$$\left(\frac{1}{2}\Delta - \frac{\partial}{\partial t}\right) v = 0 \qquad \text{on } (0, \infty) \times D,$$
$$v(0, x) = g(x), \qquad x \in D,$$
$$v(t, x) = 0, \qquad x \in \partial D.$$

In light of Theorem 1.1, if $D \subseteq \mathbb{R}^n$ is an open set, there is hope the function $u(t, x) = P_x(\tau_D(Z) > t)$ solves

$$\left(\frac{1}{8}\Delta^2 - \frac{\partial}{\partial t}\right) u(t, x) = 0 \qquad \text{on } (0, \infty) \times D,$$
$$u(0, x) = 1, \qquad x \in D,$$
$$u(t, x) = 0, \qquad x \in \partial D.$$

Since the PDE is fourth order, an additional boundary condition is needed. For this choice of $u$, it turns out to be $\frac{\partial u}{\partial n} = 0$, $x \in \partial D$, where $\frac{\partial}{\partial n}$ is the inward normal derivative. Unfortunately, this fails.

THEOREM 1.2. *Let $D = (0,1)$. Then for any $a > 0$, the function $g(t,x) = P_x(\tau_{(0,\infty)}(Z) > t)$ does not satisfy*

$$a\frac{\partial^4 g}{\partial x^4} = \frac{\partial g}{\partial t}.$$

We use another method to study the lifetime of iterated Brownian motion in an open set $D \subseteq \mathbb{R}^n$. For $n$-dimensional iterated Brownian motion $Z$, let

$$\tau_D(Z) = \inf\{t \geq 0 : Z_t \notin D\}$$

be the lifetime of $Z$ in $D$. Here and in the sequel we write $f \approx g$ and $f \lesssim g$ to mean for some positive $C_1$ and $C_2$, $C_1 \leq \frac{f}{g} \leq C_2$ and $f \leq C_1 g$, respectively. We also write $f(t) \sim g(t)$ as $t \to \infty$ to mean $\frac{f(t)}{g(t)} \to 1$ as $t \to \infty$.

Let $D \subseteq \mathbb{R}^n$ be an open cone, with vertex 0, such that $S^{n-1} \cap D$ is regular for the Laplace–Beltrami operator $L_{S^{n-1}}$ on $S^{n-1}$. We call $D$ a *generalized cone*. Then for some $p(D) > 0$, the exit time $\tau_D(B)$ of $n$-dimensional Brownian motion from $D$ satisfies

$$P_x(\tau_D(B) > t) \sim C(x) t^{-p(D)} \qquad \text{as } t \to \infty$$

([5], Corollary 1). If $n = 2$ and the angle of $D$ is $\xi \in (0, 2\pi)$, then

(1.3) $$p(D) = \frac{\pi}{2\xi};$$

see [15].

THEOREM 1.3. *Let $D \subseteq \mathbb{R}^n$ be a generalized cone. Then as $t \to \infty$*

$$P_x(\tau_D(Z) > t) \approx \begin{cases} t^{-p(D)}, & p(D) < 1, \\ t^{-1} \ln t, & p(D) = 1, \\ t^{-(p(D)+1)/2}, & p(D) > 1. \end{cases}$$

Since Brownian motion in a half-space has the same lifetime as Brownian motion in the positive reals, the next result is immediate from Theorem 1.5 and formula (1.3).

COROLLARY 1.4. *If $D \subseteq \mathbb{R}$ has the form $(-\infty, a)$ or $(b, \infty)$, then as $t \to \infty$, $P_x(\tau_D(Z) > t) \approx t^{-1/2}$.*

THEOREM 1.5. *Let $D \subseteq \mathbb{R}^n$ be bounded and open with regular boundary. If $\lambda_D$ is the principal eigenvalue of $\frac{1}{2}\Delta$ on $D$ with Dirichlet boundary conditions, then*

$$\log P_x(\tau_D(Z) > t) \sim -\tfrac{3}{2}\lambda_D^{2/3} t^{1/3} \qquad as\ t \to \infty.$$

REMARK 1.6. Compare with the Brownian case

$$P_x(\tau_D(B) > t) \sim C(x) e^{-\lambda_D t} \qquad \text{as } t \to \infty.$$

Finally, let us point out the search for probabilistic connections with higher-order PDEs has gone on for a long time (at least since the 1960s). Biharmonic functions and related boundary value problems were studied by Vanderbei [47] and Helms [25, 26]. Concerning fourth-order operators, Mądrecki [35, 36] explored connections with the 4-asymptotically stable motion and Nishioka [38, 39, 40, 41] considered the biharmonic pseudo-process—the one whose "transition density" is the fundamental solution of the parabolic biharmonic operator. Orsingher and his collaborators [1, 6, 7, 37, 42, 43, 44] studied processes governed by signed measures related to fundamental solutions of higher-order parabolic equations. Also, Hochberg and Orsingher [27] considered hyperbolic equations. Finally, Allouba [3] looked at connections with higher-order operators and Feynman–Kac-type formulas.

The article is organized as follows. The proof of Theorem 1.1 is given in Section 2. In Sections 3 and 4 the proofs of Theorems 1.3 and 1.5, respectively, are given. The last section is an appendix containing some technical results used in Sections 3 and 4. In Section 5 we prove Theorem 1.2.

**2. Proof of Theorem 1.1.** Let $f \in C_b^\infty(\mathbb{R}^n)$. As pointed out above, differentiation under the integral is not valid; hence to prove Theorem 1.1, we use distributions. Since $u$ is bounded, it yields a distribution on $\Omega = (0, \infty) \times \mathbb{R}^n$. Standard notation for this is $u \in \mathcal{D}'(\Omega)$. We show $u$ is a weak solution of $(\frac{1}{8}\Delta^2 - \frac{\partial}{\partial t})u(t,x) = -\frac{1}{2}\frac{1}{\sqrt{2\pi t}}\Delta f(x)$ on $\Omega$. Since the operator $\frac{1}{8}\Delta^2 - \frac{\partial}{\partial t}$ is hypoelliptic and $t^{-1/2}\Delta f(x) \in C^\infty(\Omega)$, by the hypoelliptic regularity theorem, $u \in C^\infty(\Omega)$ and $u$ is a classical solution to the partial differential equation. The initial condition is clearly satisfied and Theorem 1.1 will follow.

To show $u$ is a weak solution, it is enough to prove

$$(2.1) \qquad \int_\Omega u(t,x)\left[\frac{1}{8}\Delta^2 + \frac{\partial}{\partial t} + \frac{1}{2}\frac{1}{\sqrt{2\pi t}}\Delta f\right]\varphi(t,x)\,dx\,dt = 0$$

for $\varphi \in C_0^\infty(\Omega)$, the space of smooth functions with compact support in $\Omega$. We describe the case $n \geq 2$, leaving $n = 1$ to the reader. Let $\varepsilon > 0$. By

Fubini's theorem,

$$\int_\Omega u\Delta^2\varphi\,dx\,dt$$

(2.2)
$$= 2\int_\Omega \int_0^\infty \int_{B_\varepsilon(w)} p_n(y,|w-x|)f(w)p_1(t,y)\Delta^2\varphi\,dx\,dy\,dw\,dt$$
$$+ 2\int_\Omega \int_0^\infty \int_{B_\varepsilon(w)^c} p_n(y,|w-x|)f(w)p_1(t,y)\Delta^2\varphi\,dx\,dy\,dw\,dt$$
$$= I_1 + I_2, \quad \text{say,}$$

where $B_\varepsilon(w) = \{x \in \mathbb{R}^n : |w-x| < \varepsilon\}$. Now

$$\Delta_x p_n(y,|w-x|) = p_n(y,|w-x|)\left[\frac{|w-x|^2}{y^2} - \frac{n}{y}\right]$$

and, if $\frac{\partial}{\partial n_x}$ denotes the outward normal derivative on $\partial B_\varepsilon(w)^c$,

$$\frac{\partial}{\partial n_x}p_n(y,|w-x|) = \frac{\varepsilon}{y}p_n(y,\varepsilon) \qquad \text{for } x \in \partial B_\varepsilon(w).$$

Hence by Green's identity,

(2.3)
$$\int_{B_\varepsilon(w)^c} (\Delta^2\varphi)p_n(y,|w-x|)\,dx$$
$$= \int_{B_\varepsilon(w)^c} p_n(y,|w-x|)\left[\frac{|w-x|^2}{y^2} - \frac{n}{y}\right]\Delta\varphi\,dx$$
$$+ \int_{\partial B_\varepsilon(w)} \left[p_n(y,\varepsilon)\frac{\partial}{\partial n_x}\Delta\varphi - \frac{\varepsilon}{y}p_n(y,\varepsilon)\Delta\varphi\right]d\sigma(x),$$

where $\sigma(x)$ is surface measure on $\partial B_\varepsilon(w)$. Another application of Green's identity is not useful at this stage. Instead, observe

$$\frac{\partial}{\partial y}p_n(y,|w-x|) = \frac{1}{2}\Delta_w p_n(y,|w-x|)$$
$$= \frac{1}{2}\left[\frac{|w-x|^2}{y^2} - \frac{n}{y}\right]p_n(y,|w-x|).$$

Hence, (2.3) becomes

$$\int_{B_\varepsilon(w)^c} p_n(y,|w-x|)\Delta^2\varphi\,dx = 2\int_{B_\varepsilon(w)^c} \frac{\partial}{\partial y}[p_n(y,|w-x|)]\Delta\varphi\,dx$$
$$+ \int_{\partial B_\varepsilon(w)} p_n(y,\varepsilon)\left[\frac{\partial}{\partial n_x}\Delta\varphi - \frac{\varepsilon}{y}\Delta\varphi\right]d\sigma(x).$$

Now we can write $I_2$ from (2.2) as

(2.4) $$I_2 = I_3 + I_4 + I_5,$$

where
$$I_3 = 2\int_\Omega \int_0^\infty f(w)p_1(t,y)p_n(y,\varepsilon)\left[\int_{\partial B_\varepsilon(w)} \left(\frac{\partial}{\partial n_x}\Delta\varphi\right) d\sigma(x)\right] dy\, dw\, dt,$$

$$I_4 = -2\varepsilon \int_\Omega \int_0^\infty f(w)p_1(t,y)p_n(y,\varepsilon)\frac{1}{y}\left[\int_{\partial B_\varepsilon(w)} \Delta\varphi\, d\sigma(x)\right] dy\, dw\, dt,$$

$$I_5 = 4\int_\Omega \int_0^\infty \int_{B_\varepsilon(w)^c} f(w)p_1(t,y)\left[\frac{\partial}{\partial y}p_n(y,|w-x|)\right]\Delta\varphi\, dx\, dy\, dw\, dt.$$

By Fubini's theorem and integration by parts in the $dy$ integral,

$$I_5 = 4\int_\Omega \int_0^\infty \int_{B_\varepsilon(w)^c} f(w)p_1(t,y)\frac{y}{t}p_n(y,|w-x|)\Delta\varphi\, dx\, dy\, dw\, dt.$$

This has a form like $I_2$, where now $\Delta\varphi$ in $I_2$ is replaced by $\varphi$ and $p_1(t,y)$ is replaced by $\frac{2y}{t}p_1(t,y)$. Repetition of the derivation leading to (2.4) yields

(2.5) $$I_5 = I_6 + I_7 + I_8,$$

where

$$I_6 = 4\int_\Omega \int_0^\infty f(w)p_1(t,y)\frac{y}{t}p_n(y,\varepsilon)\left[\int_{\partial B_\varepsilon(w)} \frac{\partial\varphi}{\partial n_x} d\sigma(x)\right] dy\, dw\, dt,$$

$$I_7 = -4\varepsilon \int_\Omega \int_0^\infty f(w)p_1(t,y)\frac{1}{t}p_n(y,\varepsilon)\left[\int_{\partial B_\varepsilon(w)} \varphi\, d\sigma(x)\right] dy\, dw\, dt,$$

$$I_8 = 8\int_\Omega \int_0^\infty \int_{B_\varepsilon(w)^c} f(w)p_1(t,y)\frac{y}{t}\varphi\left[\frac{\partial}{\partial y}p_n(y,|w-x|)\right] dx\, dy\, dw\, dt.$$

By Fubini's theorem and integration by parts in the $dy$ integral,

$$I_8 = -8\int_\Omega \int_0^\infty \int_{B_\varepsilon(w)^c} f(w)\left[\frac{1}{t} - \frac{y^2}{t^2}\right]p_1(t,y)\varphi p_n(y,|w-x|)\, dx\, dy\, dw\, dt$$

$$= 16\int_\Omega \int_0^\infty \int_{B_\varepsilon(w)^c} f(w)\left[\frac{\partial}{\partial t}p_1(t,y)\right]\varphi p_n(y,|w-x|)\, dx\, dy\, dw\, dt$$

using that $\frac{\partial}{\partial t}p_1(t,y) = -\frac{1}{2}p_1(t,y)[\frac{1}{t} - \frac{y^2}{t^2}]$. Now we integrate by parts in the $dt$ integral to get

(2.6) $$I_8 = -16\int_\Omega \int_0^\infty \int_{B_\varepsilon(w)^c} f(w)p_1(t,y)\frac{\partial\varphi}{\partial t}p_n(y,|w-x|)\, dx\, dy\, dw\, dt.$$

Since $\varphi \in C_0^\infty(\Omega)$, there is a bounded function $h(t,w)$ with compact support in $\Omega$ such that

$$\left|\frac{\partial}{\partial n_x}\Delta\varphi(t,x)\right| \leq h(t,w), \qquad x \in \partial B_\varepsilon(w),$$

and $h$ is independent of $\varepsilon < 1$.

Then for some $M_1, M_2, M_3 > 0$,

$$|I_3| \leq C \int_\Omega \int_0^\infty h(t,w) p_1(t,y) p_n(y,\varepsilon) \sigma(\partial B_\varepsilon(w))\, dy\, dw\, dt$$

$$\leq C\varepsilon^{n-1} \int_{M_1}^{M_2} \int_{B_{M_3}(0)} \int_0^\infty p_1(t,y) p_n(y,\varepsilon)\, dy\, dw\, dt$$

$$= C\varepsilon^{n-1} \int_{M_1}^{M_2} \int_0^\infty \frac{1}{\sqrt{2\pi t}} e^{-y^2/2t} \frac{1}{(2\pi y)^{n/2}} e^{-\varepsilon^2/2y}\, dy\, dt$$

$$\leq C\varepsilon^{n-1} \int_{M_1}^{M_2} \int_0^\infty \frac{1}{\sqrt{2\pi t}} e^{-y^2/2M_2} \frac{1}{(2\pi y)^{n/2}} e^{-\varepsilon^2/2y}\, dy\, dt$$

$$= C\varepsilon^{n-1} \int_0^\infty e^{-y^2/2M_2} y^{-n/2} e^{-\varepsilon^2/2y}\, dy.$$

The charge of variables $u = \frac{\varepsilon^2}{2y}$ transforms this to

$$|I_3| \leq C\varepsilon \int_0^\infty u^{n/2-2} \exp\left(-u - \frac{\varepsilon^4}{8M_2 u^2}\right) du.$$

If $n=2$,

$$\int_0^\infty u^{-1} \exp\left(-u - \frac{\varepsilon^4}{8M_2 u^2}\right) du \leq \int_0^1 u^{-1} \exp\left(-\frac{\varepsilon^4}{8M_2 u^2}\right) du + \int_1^\infty u^{-1} e^{-u}\, du$$

$$= \frac{1}{2} \int_{\varepsilon^4/8M_2}^\infty v^{-1} e^{-v}\, dv + \int_1^\infty u^{-1} e^{-u}\, du,$$

after changing variables $v = \frac{\varepsilon^4}{8M_2} u^{-2}$ in the first integral. Using this it is easy to show $I_3 \to 0$ as $\varepsilon \to 0$. If $n \geq 3$,

$$\int_0^\infty u^{n/2-2} \exp\left(-u - \frac{\varepsilon^4}{8M_2 u^2}\right) du \leq \int_0^\infty u^{n/2-2} e^{-u}\, du < \infty.$$

Then again, $I_3 \to 0$ as $\varepsilon \to 0$. A similar argument applies to $I_6$ and $I_7$. In any event

(2.7) $$\lim_{\varepsilon \to 0} I_3 = \lim_{\varepsilon \to 0} I_6 = \lim_{\varepsilon \to 0} I_7 = 0.$$

By dominated convergence in (2.6),

(2.8) $$\lim_{\varepsilon \to 0} I_8 = -8 \int_\Omega u(t,x) \frac{\partial \varphi}{\partial t}(t,x)\, dx\, dt.$$

Also, dominated convergence gives

(2.9) $$\lim_{\varepsilon \to 0} I_1 = 0.$$

Thus (2.2) becomes

$$\int_\Omega u\Delta^2\varphi\,dx\,dt = \lim_{\varepsilon\to 0}(I_1 + I_2)$$

$$= \lim_{\varepsilon\to 0} I_2$$

(2.10)
$$= \lim_{\varepsilon\to 0}(I_3 + I_4 + I_6 + I_7 + I_8) \qquad \text{[using (2.4)–(2.5)]}$$

$$= \lim_{\varepsilon\to 0} I_4 - 8\int_\Omega u\frac{\partial\varphi}{\partial t}\,dx\,dt \qquad \text{[by (2.7)–(2.8)]}.$$

It remains to analyze $I_4$. Make the change of variables $u = \frac{\varepsilon^2}{2y}$ in the $dy$ integral to get

(2.11)
$$I_4 = -\frac{2\varepsilon^{1-n}}{\pi^{n/2}}$$
$$\times \int_\Omega \int_0^\infty \left[f(w)p_1\left(t, \frac{\varepsilon^2}{2u}\right)u^{n/2-1}e^{-u}\int_{\partial B_\varepsilon(w)}\Delta\varphi\,d\sigma(x)\right]du\,dw\,dt.$$

Since $\varphi \in C_0^\infty(\Omega)$, there is a bounded function $h_1(t,w)$ with compact support in $\Omega$ such that for $\varepsilon < 1$, the $du\,dw\,dt$ integrand (the portion square brackets) is bounded by

$$Ct^{-1/2}u^{n/2-1}e^{-u}h(t,w)\varepsilon^{n-1},$$

where $h$ and $C$ are independent of $\varepsilon < 1$. Then we can let $\varepsilon \to 0$ in (2.11) and use dominated convergence to end up with

$$\lim_{\varepsilon\to 0} I_4 = -\frac{2}{\pi^{n/2}}\int_\Omega\int_0^\infty f(w)\frac{1}{\sqrt{2\pi t}}u^{n/2-1}e^{-u}$$
$$\times \left[\lim_{\varepsilon\to 0}\varepsilon^{1-n}\int_{\partial B_\varepsilon(w)}\Delta\varphi\,d\sigma(x)\right]du\,dw\,dt$$

$$= -\frac{2}{\pi^{n/2}}\int_\Omega\int_0^\infty f(w)\frac{1}{\sqrt{2\pi t}}u^{n/2-1}e^{-u}\left[\frac{2\pi^{n/2}}{\Gamma(n/2)}\Delta\varphi(t,w)\right]du\,dw\,dt$$

$$= -4\int_\Omega f(w)\frac{1}{\sqrt{2\pi t}}\Delta\varphi(t,w)\,dw\,dt$$

$$= -4\int_\Omega \frac{1}{\sqrt{2\pi t}}\Delta f(w)\varphi(t,w)\,dw\,dt.$$

Hence, (2.10) reduces to

$$\int_\Omega u\Delta^2\varphi\,dx\,dt = -4\int_\Omega \frac{1}{\sqrt{2\pi t}}\Delta f(w)\varphi(t,w)\,dw\,dt - 8\int_\Omega u\frac{\partial\varphi}{\partial t}\,dx\,dt,$$

which is (2.1).

**3. Proof of Theorem 1.3.** If $D \subseteq \mathbb{R}$ is an open set, write

$$\tau_D^{\pm}(x) = \inf\{t \geq 0 : X_t^{\pm} + x \notin D\},$$

and if $I \subseteq \mathbb{R}$ is an open interval, write

$$\eta_I = \eta(I) = \inf\{t \geq 0 : Y_t \notin I\}.$$

By continuity of the paths,

(3.1)
$$\begin{aligned}
P_x(\tau_D > t) &\\
&= P_x(Z_s \in D \text{ for all } s \leq t) \\
&= P(x + X^+(0 \vee Y_s) \in D \text{ and } x + X^-(0 \vee (-Y_s)) \in D \text{ for all } s \leq t) \\
&= P(\tau_D^+(x) > 0 \vee Y_s \text{ and } \tau_D^-(x) > 0 \vee (-Y_s) \text{ for all } s \leq t) \\
&= P(-\tau_D^-(x) < Y_s < \tau_D^+(x) \text{ for all } s \leq t) \\
&= P(\eta(-\tau_D^-(x), \tau_D^+(x)) > t).
\end{aligned}$$

Once we show for $n$-dimensional Brownian motion $B$ and a generalized cone $D \subseteq \mathbb{R}^n$,

(3.2) $$\frac{d}{dt} P_x(\tau_D(B) \leq t) \sim c(x) t^{-p(D)-1} \quad \text{as } t \to \infty,$$

then Theorem 1.3 is an immediate consequence of (3.1) and the next theorem.

THEOREM 3.1. *Let $\xi$ be a positive random variable with density $f(t)$ such that for some positive $C$ and $p$, $f(t) \sim Ct^{-p-1}$ as $t \to \infty$. If $\xi_1$ and $\xi_2$ are independent copies of $\xi$, independent of the Brownian motion $Y$, then as $t \to \infty$,*

$$P(\eta_{(-\xi_1,\xi_2)} > t) \approx \begin{cases} t^{-p}, & p < 1, \\ t^{-1} \ln t, & p = 1, \\ t^{-(p+1)/2}, & p > 1. \end{cases}$$

We give the proof of the asymptotic (3.2) in Lemma A.3 of the Appendix.

PROOF OF THEOREM 3.1. We will abuse notation and allow $Y_0 = x$; the probability associated with this will be $P_x$. For typographical simplicity, write

$$\tau = \eta_{(0,\infty)}(Y).$$

The following distributions of $\eta_{(-u,v)}$ and $\tau$ are well known:

(3.3) $P_0(\eta_{(-u,v)} > t) = \dfrac{4}{\pi} \sum_{n=0}^{\infty} \dfrac{1}{2n+1} \exp\left(-\dfrac{(2n+1)^2\pi^2}{2(u+v)^2}t\right) \sin\dfrac{(2n+1)\pi u}{u+v},$

(3.4) $\quad P_x(\tau > t) = P_x(\eta_{(0,\infty)} > t) = \dfrac{2}{\sqrt{2\pi}} \int_0^{x/\sqrt{t}} e^{-w^2/2}\, dw$

([22], pages 340–342). By (3.4),

(3.5) $\quad\quad\quad\quad P_x(\eta_{(0,1)} > t) \leq P_x(\tau > t) \leq \dfrac{x}{\sqrt{t}} \wedge 1.$

By Lemma A.1 of the Appendix, choose $M > 0$ so large that

(3.6) $\quad P_x(\eta_{(0,1)} > t) \approx (\sin \pi x)e^{-\pi^2 t/2} \quad \text{for } t \geq M, \text{ uniformly in } x \in (0,1).$

Then choose $\delta < \frac{1}{2}$ so small that

(3.7) $\quad\quad\quad\quad \sin \pi x \approx x, \quad x \in (0, \delta].$

Choose $K > 0$ so large that

(3.8) $\quad\quad\quad\quad f(u) \approx u^{-p-1}, \quad u \geq K.$

Finally, we consider $t$ so large that

(3.9) $\quad\quad\quad\quad K < \delta\sqrt{\dfrac{t}{M}}.$

By independence of $Y, \xi_1$, and $\xi_2$, using scaling and translation invariance of Brownian motion,

$P(\eta_{(-\xi_1,\xi_2)} > t)$
$\quad = P_0(\eta_{(-\xi_1,\xi_2)} > t)$
$\quad = \displaystyle\int_0^\infty \int_0^\infty P_0(\eta_{(-u,v)} > t) f(u) f(v)\, dv\, du$
$\quad = \displaystyle\int_0^\infty \int_0^\infty P_{u/(u+v)}\left(\eta_{(0,1)} > \dfrac{t}{(u+v)^2}\right) f(u) f(v)\, dv\, du.$

Since $P_x(\eta_{(0,1)} > t) = P_{1-x}(\eta_{(0,1)} > t)$, we have

$\displaystyle\int_0^\infty \int_0^u P_{u/(u+v)}\left(\eta_{(0,1)} > \dfrac{t}{(u+v)^2}\right) f(u) f(v)\, dv\, du$

$\quad = \displaystyle\int_0^\infty \int_0^u P_{v/(u+v)}\left(\eta_{(0,1)} > \dfrac{t}{(u+v)^2}\right) f(u) f(v)\, dv\, du$

$\quad = \displaystyle\int_0^\infty \int_v^\infty P_{v/(u+v)}\left(\eta_{(0,1)} > \dfrac{t}{(u+v)^2}\right) f(u) f(v)\, du\, dv$

(reversing the order of integration)

$$= \int_0^\infty \int_u^\infty P_{u/(u+v)}\left(\eta_{(0,1)} > \frac{t}{(u+v)^2}\right) f(u)f(v)\,dv\,du$$

(relabeling). Thus,

(3.10)
$$P(\eta_{(-\xi_1,\xi_2)} > t)$$
$$= 2\int_0^\infty \int_u^\infty P_{u/(u+v)}\left(\eta_{(0,1)} > \frac{t}{(u+v)^2}\right) f(u)f(v)\,dv\,du.$$

We break up the integral into several pieces:

$$A_1 = \left\{(u,v): K \leq u \leq \delta\sqrt{\frac{t}{M}}, u+v \geq \sqrt{\frac{t}{M}}\right\},$$

$$A_2 = \left\{(u,v): 0 \leq u \leq K, u+v \geq \sqrt{\frac{t}{M}}\right\},$$

$$A_3 = \left\{(u,v): \delta\sqrt{\frac{t}{M}} \leq u \leq v, u+v \geq \sqrt{\frac{t}{M}}\right\},$$

$$A_4 = \left\{(u,v): u \geq 0, v \geq K, u \leq v, u+v \leq \sqrt{\frac{t}{M}}\right\},$$

$$A_5 = \{(u,v): u \geq 0, u \leq v \leq K\}.$$

It turns out the integral on $A_1$ is the "dominant" piece.

*The dominant piece: upper bound.* On $A_1$, $\frac{u}{u+v} \leq \frac{\delta\sqrt{t/M}}{\sqrt{t/M}} = \delta$ and so $v \geq \frac{1-\delta}{\delta}u > \frac{1-1/2}{1/2}u \geq K$. Hence, by (3.5) and (3.8),

$$\iint_{A_1} P_{u/(u+v)}\left(\eta_{(0,1)} > \frac{t}{(u+v)^2}\right) f(u)f(v)\,dv\,du$$

$$\leq \iint_{A_1} \frac{u}{\sqrt{t}} f(u)f(v)\,dv\,du$$

$$\approx t^{-1/2} \int_K^{\delta\sqrt{t/M}} \int_{\sqrt{t/M}-u}^\infty v^{-p-1} u^{-p}\,dv\,du$$

(3.11)
$$\approx t^{-1/2} \int_K^{\delta\sqrt{t/M}} u^{-p}\left(\sqrt{\frac{t}{M}} - u\right)^{-p} du$$

$$\approx t^{-1/2} \int_K^{\delta\sqrt{t/M}} u^{-p}(\sqrt{t})^{-p}\,du \qquad (t \text{ large})$$

$$= t^{-p} \int_{K/\sqrt{t}}^{\delta/\sqrt{M}} w^{-p} \, dw$$

(changing variables $u = \sqrt{t} w$).

*The dominant piece: lower bound.* Notice on $A_1$, $\frac{u}{u+v} \leq \frac{\delta\sqrt{t/M}}{\sqrt{t/M}} = \delta$ and $\frac{t}{(u+v)^2} \leq M$. Hence,

$$\iint_{A_1} P_{u/(u+v)}\left(\eta_{(0,1)} > \frac{t}{(u+v)^2}\right) f(u) f(v) \, dv \, du$$

$$\geq \iint_{A_1} P_{u/(u+v)}(\eta_{(0,1)} > M) f(u) \, du \, f(v) \, dv$$

$$\approx \int_K^{\delta\sqrt{t/M}} \int_{\sqrt{t/M}-u}^{\infty} \frac{u}{u+v} v^{-p-1} u^{-p-1} \, dv \, du \qquad [\text{by (3.6)--(3.8)}]$$

$$= t^{-p} \int_{K/\sqrt{t}}^{\delta/\sqrt{M}} \int_{1/\sqrt{M}}^{\infty} \frac{w}{z} (z-w)^{-p-1} w^{-p-1} \, dz \, dw$$

[changing variables $z = (u+v)/\sqrt{t}, w = u/\sqrt{t}$]

$$\geq t^{-p} \int_{K/\sqrt{t}}^{\delta/\sqrt{M}} \int_{1/\sqrt{M}}^{\infty} \frac{w}{z} \left(z - \frac{K}{\sqrt{t}}\right)^{-p-1} w^{-p-1} \, dz \, dw$$

$$= t^{-p} \int_{1/\sqrt{M}}^{\infty} \frac{1}{z} \left(z - \frac{K}{\sqrt{t}}\right)^{-p-1} dz \int_{K/\sqrt{t}}^{\delta/\sqrt{M}} w^{-p} \, dw$$

$$\approx t^{-p} \int_{1/\sqrt{M}}^{\infty} z^{-p-2} \, dz \int_{K/\sqrt{t}}^{\delta/\sqrt{M}} w^{-p} \, dw \qquad \text{for large } t.$$

Combined with (3.11), we end up with

(3.12)
$$\iint_{A_1} P_{u/(u+v)}\left(\eta_{(0,1)} > \frac{t}{(u+v)^2}\right) f(u) f(v) \, dv \, du$$
$$\approx t^{-p} \int_{K/\sqrt{t}}^{\delta/\sqrt{M}} w^{-p} \, dw.$$

*The remaining terms: upper bounds.* On the set $A_2$, $v \geq \sqrt{\frac{t}{M}} - K > K$ if $t$ is large enough. Hence, by (3.5) and (3.8),

$$\iint_{A_2} \lesssim \int_0^K \int_{\sqrt{t/M}-u}^{\infty} \frac{u}{\sqrt{t}} v^{-p-1} \, dv \, f(u) \, du$$

$$\approx \int_0^K t^{-1/2} u \left(\sqrt{\frac{t}{M}} - u\right)^{-p} f(u) \, du$$

(3.13)
$$\leq t^{-1/2}\left(\sqrt{\frac{t}{M}} - K\right)^{-p} \int_0^K u f(u)\, du$$
$$\approx t^{-(p+1)/2} \quad \text{for } t \text{ large.}$$

On the set $A_3$, if $t$ is large, then $v \geq u \geq K$ so by (3.5) and (3.8),

$$\iint_{A_3} \leq \iint_{A_3} 1 \cdot f(u)f(v)\, dv\, du \approx \iint_{A_3} u^{-p-1} v^{-p-1}\, dv\, du$$

$$= \int_{\delta\sqrt{t/M}}^{(\sqrt{t/M})/2} \int_{\sqrt{t/M}-u}^{\infty} v^{-p-1} u^{-p-1}\, dv\, du$$

(3.14)
$$+ \int_{(\sqrt{t/M})/2}^{\infty} \int_u^{\infty} v^{-p-1} u^{-p-1}\, dv\, du$$

$$\leq \int_{\delta\sqrt{t/M}}^{(\sqrt{t/M})/2} \frac{1}{p}\left(\frac{1}{2}\sqrt{\frac{t}{M}}\right)^{-p} u^{-p-1}\, du + \int_{(\sqrt{t/M})/2}^{\infty} \frac{1}{p} u^{-2p-1}\, du$$

$$\approx t^{-p} + t^{-p}$$
$$\approx t^{-p}.$$

On $A_4$, since $\frac{t}{(u+v)^2} \geq M$, by (3.6) and that $\frac{\sin \pi x}{x}$ is bounded for all $x$,

$$\iint_{A_4} P_{u/(u+v)}\left(\eta_{(0,1)} > \frac{t}{(u+v)^2}\right) f(u)f(v)\, dv\, du$$

$$\lesssim \iint_{A_4} \frac{u}{u+v} e^{-\pi^2 t/2(u+v)^2} f(u)f(v)\, dv\, du$$

$$\approx \int_0^K \int_K^{\sqrt{t/M}-u} \frac{u}{u+v} e^{-\pi^2 t/2(u+v)^2} f(u) v^{-p-1}\, dv\, du$$

$$+ \int_K^{\sqrt{t/M}/2} \int_u^{\sqrt{t/M}-u} \frac{u}{u+v} e^{-\pi^2 t/2(u+v)^2} u^{-p-1} v^{-p-1}\, dv\, du$$

[using (3.8)]. Changing variables $z = \frac{u+v}{\sqrt{t}}$, $w = u$ in the first integral and $z = \frac{u+v}{\sqrt{t}}$, $w = \frac{u}{\sqrt{t}}$ in the second to get

$$\int_0^K \int_{(w+K)/\sqrt{t}}^{1/\sqrt{M}} \frac{w}{z\sqrt{t}} e^{-\pi^2/2z^2} f(w)(z\sqrt{t} - w)^{-p-1}\sqrt{t}\, dz\, dw$$

$$+ \int_{K/\sqrt{t}}^{1/(2\sqrt{M})} \int_{2w}^{1/\sqrt{M}} \frac{w}{z} e^{-\pi^2/2z^2} (w\sqrt{t})^{-p-1}(z\sqrt{t} - w\sqrt{t})^{-p-1} t\, dz\, dw$$

$$= t^{-(p+1)/2} \int_0^K \int_{(w+K)/\sqrt{t}}^{1/\sqrt{M}} \frac{w}{z} e^{-\pi^2/2z^2} f(w) z^{-p-1}\left(1 - \frac{w}{z\sqrt{t}}\right)^{-p-1} dz\, dw$$

$$+ t^{-p} \int_{K/\sqrt{t}}^{1/(2\sqrt{M})} \int_{2w}^{1/\sqrt{M}} \frac{1}{z} e^{-\pi^2/2z^2} w^{-p} z^{-p-1} \left(1 - \frac{w}{z}\right)^{-p-1} dz\, dw.$$

In the first term, since $0 \leq w \leq K$ and $\frac{w+K}{\sqrt{t}} \leq z$, we have

$$\frac{w}{z\sqrt{t}} \leq \frac{w}{w+K} \leq \frac{1}{2};$$

hence, $(1 - \frac{w}{z\sqrt{t}})^{-p-1} < (\frac{1}{2})^{-p-1}$. In the second term, since $2w < z$, we have

$$\frac{w}{z} < \frac{1}{2};$$

hence, $(1 - \frac{w}{z})^{-p-1} < (\frac{1}{2})^{-p-1}$.

It follows that

$$\iint_{A_4} \lesssim t^{-(p+1)/2} \int_0^K \int_{(w+K)/\sqrt{t}}^{1/\sqrt{M}} wz^{-p-2} e^{-\pi^2/2z^2} f(w)\, dz\, dw$$

$$+ t^{-p} \int_{K/\sqrt{t}}^{1/(2\sqrt{M})} \int_{2w}^{1/\sqrt{M}} z^{-p-2} w^{-p} e^{-\pi/2z^2} dz\, dw$$

$$\lesssim t^{-(p+1)/2} \int_0^K \int_0^{1/\sqrt{M}} wz^{-p-2} e^{-\pi^2/2z^2} f(w)\, dz\, dw$$

(3.15)
$$+ t^{-p} \left(\int_0^{1/\sqrt{M}} z^{-p-2} e^{-\pi^2/2z^2} dz\right) \int_{K/\sqrt{t}}^{1/(2\sqrt{M})} w^{-p}\, dw$$

$$\approx t^{-(p+1)/2} + t^{-p} \int_{K/\sqrt{t}}^{1/(2\sqrt{M})} w^{-p}\, dw$$

$$\approx t^{-(p+1)/2} + \iint_{A_1} + t^{-p},$$

using (3.12).

Finally, on $A_5$, since $\frac{t}{(u+v)^2} \geq M$, by (3.6),

(3.16)
$$\iint_{A_5} \leq \int_0^K \int_0^K e^{-\pi^2 t/2(u+v)^2} f(u)f(v)\, dv\, du \leq Ce^{-\pi^2 t/8K^2}$$

$$\leq t^{-p} \quad \text{for } t \text{ large}.$$

Using (3.10) and (3.13)–(3.16), for $t$ large,

$$\left(\iint_{A_1}\right) \lesssim P(\eta_{(-\xi_1,\xi_2)} > t)$$

$$\lesssim \left[\left(\iint_{A_1}\right) + t^{-(p+1)/2} + t^{-p}\right]$$

$$+ \left( t^{-(p+1)/2} + \left( \iint_{A_1} \right) + t^{-p} \right) + t^{-p} \bigg].$$

By (3.12), for $t$ large,

$$\left( \iint_{A_1} \right) \approx \begin{cases} t^{-p}, & p < 1, \\ t^{-1} \ln t, & p = 1, \\ t^{-(p+1)/2}, & p > 1. \end{cases}$$

Hence, the desired conclusion follows. □

**4. Exponential decay.** In this section we prove Theorem 1.5. First we state a couple of lemmas.

LEMMA 4.1. *Let $\xi$ be a positive random variable such that for some $c > 0$, $-\log P(\xi > t) \sim ct$ as $t \to \infty$. Then for independent copies $\xi_1$ and $\xi_2$ of $\xi$, $-\log P(\xi_1 + \xi_2 > t) \sim ct$ as $t \to \infty$.*

PROOF. By independence,

$$P(\xi_1 + \xi_2 > t) = 1 - \int_0^t P(\xi \leq t - y) \, d_y P(\xi \leq y)$$

$$= 1 - \int_0^t [1 - P(\xi > t - y)] \, d_y P(\xi \leq y)$$

$$= P(\xi > t) + \int_0^t P(\xi > t - y) \, d_y P(\xi \leq y)$$

$$\geq P(\xi > t).$$

For an upper bound, note for any $\theta < c$, $E[e^{\theta \xi}] < \infty$. Then by the Markov–Chebyshev inequality and independence,

$$P(\xi_1 + \xi_2 > t) \leq e^{-\theta t} E[e^{\theta(\xi_1 + \xi_2)}]$$

$$= e^{-\theta t} [E[e^{\theta \xi}]]^2.$$

Thus,

$$\theta t - 2 \log E[e^{\theta \xi}] \leq -\log P(\xi_1 + \xi_2 > t) \leq -\log P(\xi > t).$$

Divide by $ct$, let $t \to \infty$, then let $\theta \uparrow c$ to get the desired conclusion. □

We use Lemma 4.1 to derive asymptotics of the Laplace transform of $(\xi_1 + \xi_2)^{-2}$. This is a simple application of de Bruijn's Tauberian theorem ([9], page 254). We state a special case of it for the convenience of the reader.

DE BRUIJN'S TAUBERIAN THEOREM (Special case). *Let $X$ be a positive random variable such that for some positive $B$ and $p$,*

$$-\log P(X \leq x) \sim Bx^{-p} \qquad \text{as } x \to 0.$$

*Then*

$$-\log Ee^{-\lambda X} \sim (p+1)B^{1/(p+1)}p^{-p/(p+1)}\lambda^{p/(p+1)} \qquad \text{as } \lambda \to \infty.$$

LEMMA 4.2. *With $\xi_1$ and $\xi_2$ as in Lemma 4.1,*

$$-\log E(\exp(-\lambda/(\xi_1+\xi_2)^2)) \sim 3(c/2)^{2/3}\lambda^{1/3} \qquad \text{as } \lambda \to \infty.$$

PROOF. By Lemma 4.2,

$$-\log P\left(\frac{1}{(\xi_1+\xi_2)^2} \leq x\right) = -\log P(\xi_1+\xi_2 > x^{-1/2})$$

$$\sim cx^{-1/2} \qquad \text{as } x \to 0.$$

The claim follows from de Bruijn's theorem. □

We will also need the following application of de Bruijn's theorem.

LEMMA 4.3. *Let $X$ be a positive random variable with density $f(u) = \gamma u^{-2}e^{-\alpha/\sqrt{u}}$ for $u > 0$. Here $\alpha > 0$ and $\gamma$ is chosen to give total mass 1. Then $-\log Ee^{-\lambda X} \sim 3\alpha^{2/3}2^{-2/3}\lambda^{1/3}$ as $\lambda \to \infty$.*

PROOF. We have, changing variables $u = v^{-2}$,

$$P(X \leq x) = \int_0^x \gamma u^{-2}e^{-\alpha/\sqrt{u}}\,du$$

$$= \int_{x^{-1/2}}^\infty 2\gamma v e^{-\alpha v}\,dv$$

$$= \frac{2\gamma}{\alpha^2}[\alpha x^{-1/2}+1]e^{-\alpha x^{-1/2}}.$$

Hence,

$$-\log P(X \leq x) \sim \alpha x^{-1/2} \qquad \text{as } x \to 0.$$

The desired conclusion follows from de Bruijn's theorem. □

The main result of this section is the next theorem. Before proving it, we show how it yields Theorem 1.5.

THEOREM 4.4. *Let $\xi > 0$ be a positive random variable whose density $f$ satisfies $-\log f(t) \sim ct$ as $t \to \infty$. If $\xi_1$ and $\xi_2$ are independent copies of $\xi$, independent of the Brownian motion $Y$, then*

$$-\log P(\eta_{(-\xi_1, \xi_2)} > t) \sim \tfrac{3}{2} c^{2/3} t^{1/3} \qquad as\ t \to \infty.$$

In Lemma A.4 of the Appendix we prove for a bounded open set $D \subseteq \mathbb{R}^n$ with regular boundary,

$$\log\left(\frac{d}{dt} P_x(\tau_D(B) \le t)\right) \sim -\lambda_D t \qquad as\ t \to \infty,$$

where $\lambda_D$ is the principal eigenvalue of $\tfrac{1}{2}\Delta$ on $D$ with Dirichlet boundary condition. Then by (3.1) and Theorem 4.4, Theorem 1.5 follows.

PROOF OF THEOREM 4.4. First we give a lower bound. Let $\varepsilon > 0$. Recall we use the notation $g \gtrsim f$ to mean for some constant $c_1 > 0$, $g \ge c_1 f$. Choose $M$ and $\delta$ as in (3.6) and (3.7). By hypotheses, choose $K > 0$ such that

(4.1) $$e^{-c(1+\varepsilon)u} \le f(u) \le e^{-c(1-\varepsilon)u} \qquad \text{for } u \ge K.$$

For $t$ large, let

$$A = \left\{(u,v) : K \le u \le \delta\sqrt{\frac{t}{M}},\ \frac{1-\delta}{\delta} u \le v \le \sqrt{\frac{t}{M}} - u\right\}.$$

Then by (3.10),

$$P(\eta_{(-\xi_1, \xi_2)} > t) \ge \iint_A P_{u/(u+v)}\left(\eta_{(0,1)} > \frac{t}{(u+v)^2}\right) f(u) f(v)\, dv\, du.$$

On the set $A$, since $\delta < \tfrac{1}{2}$, we have $v \ge (\tfrac{1}{\delta} - 1) u > u > K$ and $u + v > \tfrac{u}{\delta}$. Thus, $\frac{u}{u+v} \le \delta$ and so by (3.6), (3.7) and (4.1),

$$P(\eta_{(-\xi_1,\xi_2)} > t)$$
$$\gtrsim \int_K^{\delta\sqrt{t/M}} \int_{(1-\delta)u/\delta}^{\sqrt{t/M}-u} \frac{u}{u+v} e^{-\pi^2 t/2(u+v)^2} e^{-c(1+\varepsilon)u} e^{-c(1+\varepsilon)v}\, dv\, du$$

(changing variables $x = u+v, z = u$)

$$\approx \int_K^{\delta\sqrt{t/M}} \int_{z/\delta}^{\sqrt{t/M}} \frac{z}{x} e^{-\pi^2 t/2x^2} e^{-c(1+\varepsilon)x}\, dx\, dz$$

(reversing the order of integration)

$$= \int_{K/\delta}^{\sqrt{t/M}} \int_K^{\delta x} \frac{z}{x} e^{-\pi^2 t/2x^2} e^{-c(1+\varepsilon)x} \, dz \, dx$$

$$\approx \int_{K/\delta}^{\sqrt{t/M}} \frac{1}{x} e^{-\pi^2 t/2x^2} e^{-c(1+\varepsilon)x} [\delta^2 x^2 - K^2] \, dx$$

$$\geq \int_{2K/\delta}^{\sqrt{t/M}} \frac{1}{x} e^{-\pi^2 t/2x^2} e^{-c(1+\varepsilon)x} (\delta^2 x^2 - K^2) \, dx \qquad \text{(for } t \text{ large)}$$

$$\geq \int_{2K/\delta}^{\sqrt{t/M}} \frac{1}{x} e^{-\pi^2 t/2x^2} e^{-c(1+\varepsilon)x} \left(\delta^2 x^2 - \frac{\delta^2 x^2}{4}\right) dx$$

(changing variables $u = x^{-2}$)

$$\approx \int_{M/t}^{\delta^2/4K^2} u^{-2} e^{-c(1+\varepsilon)/\sqrt{u}} e^{-\pi^2 tu/2} \, du.$$

Thus we have for $t$ large,

$$(4.2) \qquad P(\eta_{(-\xi_1, \xi_2)} > t) \gtrsim \int_{M/t}^{\delta^2/4K^2} u^{-2} e^{-c(1+\varepsilon)/\sqrt{u}} e^{-\pi^2 tu/2} \, du.$$

Now

$$(4.3) \qquad \begin{aligned} &\int_{\delta^2/4K^2}^{\infty} u^{-2} e^{-c(1+\varepsilon)/\sqrt{u}} e^{-\pi^2 tu/2} \, du \\ &\leq e^{-\pi \delta^2 t/8K^2} \int_{\delta^2/4K^2}^{\infty} u^{-2} e^{-c(1+\varepsilon)/\sqrt{u}} \, du \lesssim e^{-c_1 t} \end{aligned}$$

and changing variables $u = v^{-2}$,

$$(4.4) \qquad \begin{aligned} &\int_0^{M/t} u^{-2} e^{-c(1+\varepsilon)/\sqrt{u}} e^{-\pi^2 tu/2} \, du \\ &\leq \int_0^{M/t} u^{-2} e^{-c(1+\varepsilon)/\sqrt{u}} \, du \\ &= 2 \int_{\sqrt{t/M}}^{\infty} v e^{-c(1+\varepsilon)v} \, dv \lesssim \sqrt{t} e^{-c(1+\varepsilon)\sqrt{t/M}}. \end{aligned}$$

By Lemma 4.3, with $\alpha = c(1+\varepsilon)$ and $\lambda = \frac{\pi^2 t}{2}$,

$$(4.5) \qquad -\log \int_0^{\infty} u^{-2} e^{-c(1+\varepsilon)/\sqrt{u}} e^{-\pi^2 tu/2} \, du \sim \tfrac{3}{2} c^{2/3} \pi^{2/3} (1+\varepsilon)^{2/3} t^{1/3}$$

$$\text{as } t \to \infty.$$

By (4.3)–(4.5), we can rewrite (4.2) for large $t$ as

(4.6) $\quad P(\eta_{(-\xi_1,\xi_2)} > t) \gtrsim \exp(-\tfrac{3}{2}c^{2/3}\pi^{2/3}(1+\varepsilon)^{5/3}t^{1/3}).$

Now we give an upper bound. By (3.10) and (3.6),

$$P(\eta_{(-\xi_1,\xi_2)} > t) = 2\left[\iint_{u+v \leq \sqrt{t/M}} + \iint_{u+v \geq \sqrt{t/M}}\right]$$

$$\times \left(P_{u/(u+v)}\left(\eta_{(0,1)} > \frac{t}{(u+v)^2}\right)f(u)f(v)\,dv\,du\right)$$

(4.7) $\quad \lesssim \iint_{u+v \leq \sqrt{t/M}} e^{-\pi^2 t/2(u+v)^2} f(u)f(v)\,dv\,du$

$$+ \iint_{u+v \geq \sqrt{t/M}} f(u)f(v)\,dv\,du$$

$$\leq E\left[\exp\left(-\frac{\pi^2 t}{2(\xi_1+\xi_2)^2}\right)\right] + P(\xi_1+\xi_2 \geq \sqrt{t/M}).$$

Since $P(\xi > t) = \int_t^\infty f(u)\,du$, by our hypotheses on the density $f$ of $\xi$,

$$\log P(\xi > t) \sim -ct \quad \text{as } t \to \infty.$$

By Lemmas 4.1 and 4.2, for $t$ large (4.7) becomes

$$P(\eta_{(-\xi_1,\xi_2)} > t) \lesssim \exp(-\tfrac{3}{2}c^{2/3}\pi^{2/3}(1-\varepsilon)t^{1/3})$$

$$+ \exp(-c(1-\varepsilon)\sqrt{t/M})$$

$$\lesssim \exp(-\tfrac{3}{2}c^{2/3}\pi^{2/3}(1-\varepsilon)t^{1/3}).$$

Combined with (4.6), this yields

$$-\tfrac{3}{2}c^{2/3}\pi^{2/3}(1+\varepsilon)^{5/3} \leq \liminf_{t\to\infty} t^{-1/3}\log P(\eta_{(-\xi_1,\xi_2)} > t)$$

$$\leq \limsup_{t\to\infty} t^{-1/3}\log P(\eta_{(-\xi_1,\xi_2)} > t)$$

$$\leq -\tfrac{3}{2}c^{2/3}\pi^{2/3}(1-\varepsilon).$$

Let $\varepsilon \to 0$ to get the desired conclusion. □

**5. Proof of Theorem 1.2.** In one dimension, with the boundary conditions $0 = f(0) = f(1) = f'(0) = f'(1)$, the biharmonic operator $\frac{d^4}{dx^4}$ is known to have a complete orthonormal set of eigenfunctions $\{\varphi_n : n \geq 1\}$ in $L^2(0,1)$, with corresponding eigenvalues $\{\lambda_n : n \geq 1\}$. In fact, the eigenvalues are of the form $\lambda = \alpha^4$, where $\alpha$ is a positive solution of

$$\cos(\alpha)\cosh(\alpha) = 1,$$

and the corresponding eigenfunction is

$$
\begin{aligned}
(5.1) \quad c(\alpha)[&(\sin(\alpha) + \sinh(\alpha))(\sin(\alpha x) - \sinh(\alpha x)) \\
&+ (\cos(\alpha) - \cosh(\alpha))(\cos(\alpha x) - \cosh(\alpha x))],
\end{aligned}
$$

with $c(\alpha)$ chosen to make the $L^2(0,1)$ norm one. Moreover,

$$(5.2) \qquad \lambda_n \sim \left[\frac{(2n+1)\pi}{2}\right]^4 \qquad \text{as } n \to \infty.$$

These facts can be found in [19] on pages 113–116.

LEMMA 5.0. *The eigenfunctions are uniformly bounded:*

$$\sup_{n,x} |\varphi_n(x)| < \infty.$$

PROOF. Simple computation reveals the constant $c(\alpha)$ in (5.1) satisfies $c(\alpha) \sim 2e^{-\alpha}$ as $\alpha \to \infty$. Expansion of (5.1) reveals the only trouble comes from the

$$-\sinh(\alpha)\sinh(\alpha x) + \cosh(\alpha)\cosh(\alpha x)$$

terms. However, these combine to give $\cosh(\alpha(1-x))$ and clearly $c(\alpha)\cosh(\alpha(1-x))$ is bounded for $x \in (0,1)$. □

From (3.1), since $P(\tau^\pm_{(0,1)}(x) > t) = P_x(\eta_{(0,1)} > t)$, we have

$$(5.3) \quad P_x(\tau_{(0,1)}(Z) > t) = \int_0^\infty \int_0^\infty P_0(\eta_{(-u,v)} > t) f(x,u) f(x,u) \, du \, dv,$$

where $f(x,u) = -\frac{d}{du} P_x(\eta_{(0,1)} > u)$ is the density of $\eta_{(0,1)}$. This density has two expansions—one good for large $u$, the other good for small $u$:

$$(5.4) \quad f(x,u) = 2\pi \sum_{n=0}^\infty (2n+1) \exp\left(-\frac{(2n+1)^2 \pi^2}{2} u\right) \sin((2n+1)\pi x),$$

$$
\begin{aligned}
(5.5) \quad f(x,u) = &(2\pi)^{-1/2} u^{-3/2} \\
&\times \sum_{k=-\infty}^\infty [(x+2k)e^{-(x+2k)^2/2u} + (1-x+2k)e^{-(1-x+2k)^2/2u}]
\end{aligned}
$$

[see [22], page 342, and [10], page 172, (3.0.2)]. These expansions can be used to prove the following lemmas. We omit the details.

LEMMA 5.1. *For $g(t,x) = P_x(\tau_{(0,1)}(Z) > t)$, the derivatives $\frac{\partial^n g}{\partial x^n}$ are bounded for $x$ in compact subsets of $(0,1)$, $t > 0$ and $n \leq 4$.*

LEMMA 5.2. *The function $g(t,x) = P_x(\tau_{(0,1)}(Z) > t)$ satisfies the boundary conditions $\frac{\partial^n g}{\partial x^n} = 0$ at $x = 0, 1$ for $n \leq 1$.*

Now we prove Theorem 1.2. To get a contradiction, assume

$$g(t,x) = P_x(\tau_{(0,1)}(Z) > t)$$

satisfies the equation

$$\frac{\partial g}{\partial t} = a \frac{\partial^4 g}{\partial x^4}.$$

Taking the Laplace transform in $t$ and using Lemma 5.1 to interchange integration and differentiation, this becomes

$$\lambda \hat{g}(\lambda, x) - 1 = a \frac{\partial^4}{\partial x^4} \hat{g}(\lambda, x), \tag{5.6}$$

where

$$\hat{g}(\lambda, x) = \int_0^\infty e^{-\lambda t} g(t, x)\, dt.$$

Since $g$ is bounded, so is $\hat{g}$; hence for each $\lambda > 0$, we can expand with respect to the orthonormal basis $\{\varphi_n : n \geq 1\}$ of $L^2(0,1)$,

$$\hat{g}(\lambda, x) = \sum_{n=1}^\infty c_n(\lambda) \varphi_n(x) \quad \text{in } L^2(0,1),$$

where

$$c_n(\lambda) = \int_0^1 \hat{g}(\lambda, x) \varphi_n(x)\, dx.$$

To determine $c_n(\lambda)$, note by (5.6)

$$\lambda c_n(\lambda) - \int_0^1 \varphi_n(y)\, dy = a \int_0^1 \left[\frac{\partial^4}{\partial x^4} \hat{g}(\lambda, x)\right] \varphi_n(x)\, dx.$$

By Lemmas 5.1 and 5.2, $\frac{\partial^n}{\partial x^n} \hat{g}(\lambda, x) = 0$ for $x = 0, 1$ and $n \leq 1$. Since $\varphi_n$ satisfies the same boundary conditions, we can integrate by parts four times to get

$$\lambda c_n(\lambda) - \int_0^1 \varphi_n(y)\, dy = a \int_0^1 \hat{g}(\lambda, x) \varphi_n^{(4)}(x)\, dx$$

$$= -a \int_0^1 \hat{g}(\lambda, x) \lambda_n \varphi_n(x)\, dx$$

$$= -a \lambda_n c_n(\lambda).$$

Thus,
$$c_n(\lambda) = \frac{\int_0^1 \varphi_n(y)\,dy}{\lambda + a\lambda_n}.$$

Consequently,
$$\hat{g}(\lambda, x) = \sum_{n=1}^{\infty} \frac{\int_0^1 \varphi_n(y)\,dy}{\lambda + a\lambda_n} \varphi_n(x) \qquad \text{in } L^2(0,1).$$

By (5.2) and Lemma 5.0, the convergence is also pointwise. Then by the continuity theorem for Laplace transforms, for $0 < a < b$,

$$(5.7) \qquad \sum_{n=1}^{\infty} \int_a^b \left[ \int_0^1 \varphi_n(y)\,dy \right] \varphi_n(x) e^{-a\lambda_n s}\,ds = \int_a^b g(s, x)\,ds.$$

For $t$ large, say $t \geq T$,
$$\sum_{n=1}^{\infty} \left[ \int_0^1 \varphi_n(y)\,dy \right] \varphi_n(x) e^{-a\lambda_n t} \approx \left[ \int_0^1 \varphi_1(y)\,dy \right] \varphi_1(x) e^{-a\lambda_1 t}$$
$$\approx C(x) e^{-a\lambda_1 t}$$

[by Lemma 5.0 and (5.2)] and
$$\exp(-\tfrac{3}{2} \cdot \tfrac{3}{2} \lambda_{(0,1)}^{2/3} t^{1/3}) \leq g(t, x)$$

(by Theorem 1.5). Hence for $t \geq T$, (5.7) yields
$$\int_t^{2t} \exp(-c_1 s^{1/3})\,ds \lesssim \int_t^{2t} e^{-a\lambda_1 s}\,ds.$$

This implies
$$t \exp(-c_1 2^{1/3} t^{1/3}) \leq c_2 e^{-a\lambda_1 t} [1 - e^{-a\lambda_1 t}], \qquad t \geq T.$$

Take the natural logarithm of both sides, divide by $t^{1/3}$ and let $t \to \infty$ to get
$$-c_1 2^{1/3} \leq -\infty,$$
a contradiction.

## APPENDIX

We collect some technical results used above.

LEMMA A.1. *As $t \to \infty$,*
$$P_x(\eta_{(0,1)} > t) \sim \frac{4}{\pi} e^{-\pi^2 t/2} \sin \pi x \qquad \text{uniformly for } x \in (0,1).$$

Proof. In light of the eigenfunction expansion

$$P_x(\eta_{(0,1)} > t) = \frac{4}{\pi} \sum_{n=0}^{\infty} \frac{1}{(2n+1)} e^{-2(n+1)^2 \pi^2 t/2} \sin(2n+1)\pi x$$

([22], page 342), it suffices to show for some $c > 0$, independent of $n$,

$$\frac{|\sin(2n+1)\pi x|}{\sin \pi x} \leq c(2n+1)^2, \qquad x \in (0,1).$$

Since

$$\frac{|\sin(2n+1)\pi(1-x)|}{\sin \pi(1-x)} = \frac{|\sin(2n+1)\pi x|}{\sin \pi x},$$

we can restrict attention to $x \in (0, \frac{1}{2}]$. There exists $c > 0$ such that $\sin \pi x \geq cx$, $x \in (0, \frac{1}{2}]$ and so

$$\frac{|\sin(2n+1)\pi x|}{\sin \pi x} \leq \frac{(2n+1)\pi x}{cx} = \frac{\pi}{c}(2n+1), \qquad x \in \left(0, \frac{1}{2}\right]. \qquad \square$$

To prove (3.2), we collect some facts. If $D \subseteq \mathbb{R}^n$ is a generalized cone, there is a complete set of orthonormal eigenfunctions $m_j$ with corresponding eigenvalues $0 < \lambda_1 \leq \lambda_2 \leq \cdots$ satisfying

$$L_{S^{n-1}} m_j = -\lambda_j m_j \qquad \text{on } D \cap S^{n-1},$$
$$m_j = 0 \qquad \text{on } S^{n-1} \cap \partial D$$

and for $\sigma(d\theta)$ surface measure on $S^{n-1}$,

$$\int_{D \cap S^{n-1}} m_j^2(\theta) \sigma(d\theta) = 1$$

[16]. For notational simplicity, set

(A.1) $$\gamma_j = \sqrt{\lambda_j + \left(\frac{n}{2} - 1\right)^2}.$$

The confluent hypergeometric function (for $b > 0$) is

$$_1F_1(a,b,z) = 1 + \frac{a}{b}\frac{z}{1!} + \frac{a(a+1)}{b(b+1)}\frac{z^2}{2!} + \cdots;$$

Bañuelos and Smits [5] have shown there is a series expansion

(A.2) $$P_x(\tau_D(B) > t) = \sum_{j=1}^{\infty} B_j \left(\frac{|x|^2}{2t}\right)^{a_j/2} {}_1F_1\left(\frac{a_j}{2}, a_j + \frac{n}{2}, -\frac{|x|^2}{2t}\right) m_j\left(\frac{x}{|x|}\right),$$

which converges uniformly for $(x,t) \in K \times [T, \infty)$, where $K \subseteq D$ is compact, $T > 0$,

(A.3) $$a_j = \gamma_j - \left(\frac{n}{2} - 1\right)$$

and

(A.4) $$B_j = \frac{\Gamma((a_j + n)/2)}{\Gamma(a_j + n/2)} \int_{D \cap S^{n-1}} m_j(\theta)\sigma(d\theta).$$

They also show for each $x \in D$,

$$P_x(\tau_D(B) > t) \sim B_1 m_1\left(\frac{x}{|x|}\right)\left(\frac{|x|^2}{2}\right)^{a_1/2} t^{-a_1/2} \quad \text{as } t \to \infty.$$

Hence, the $p(D)$ described in the Introduction and appearing in the statement of Theorem 1.5 is $a_1/2$.

LEMMA A.2. *Let $k > 0$ be an integer. If $a > k$ and $\gamma - \alpha - 1 \geq 0$, then*

$$0 \leq x^\alpha {}_1F_1(\alpha, \gamma, -x) \leq x^k \frac{\Gamma(\gamma)\Gamma(\alpha - k)}{\Gamma(\alpha)\Gamma(\gamma - \alpha)}, \quad x \geq 0.$$

PROOF. By formula 9.211.2 on page 1058 of [24],

$$x^\alpha \frac{\Gamma(\alpha)\Gamma(\gamma - \alpha)}{\Gamma(\gamma)} {}_1F_1(\alpha, \gamma, -x)$$

$$= x^\alpha(-x)^{1-\gamma} \int_0^{-x} e^t t^{\alpha-1}(-x-t)^{\gamma-\alpha-1}\, dt$$

$$= x^{\alpha+1-\gamma} \int_0^x e^{-u} u^{\alpha-1} x^{\gamma-\alpha-1}\left(1 - \frac{u}{x}\right)^{\gamma-\alpha-1} du$$

$$= \int_0^x e^{-u} u^{\alpha-1}\left(1 - \frac{u}{x}\right)^{\gamma-\alpha-1} du$$

$$= \int_0^x e^{-u} u^{\alpha-1-k} u^k \left(1 - \frac{u}{x}\right)^{\gamma-\alpha-1} du$$

$$\leq x^k \int_0^x e^{-u} u^{\alpha-1-k}\, du \quad (\text{since } \gamma - \alpha - 1 \geq 0)$$

$$\leq x^k \Gamma(\alpha - k).$$

It is clear from the nonnegativity of the integral in the third equality above that $x^\alpha {}_1F_1(\alpha, \gamma, -x) \geq 0$. □

From formula 9.213 on page 1058 in [24],

(A.5) $$\frac{d}{dx} {}_1F_1(\alpha, \gamma, x) = \frac{\alpha}{\gamma} {}_1F_1(\alpha + 1, \gamma + 1, x)$$

and by Theorem 8, page 102, in [16],

(A.6) $$\sup_{D \cap S^{n-1}} |m_j|^2 \leq C(n) \lambda_j^{(n-1)/2}.$$

Also, the proof of Theorem 4.1 part (iv) in [20] carries through in our context to yield

(A.7) $$\sum_{j=1}^{\infty} \lambda_j^{-2} < \infty.$$

LEMMA A.3.  *The asymptotic* (3.2) *is valid.*

PROOF.  Let $f(t)$ be the density of $P_x(\tau_D(B) \leq t)$. Then

$$f(t) = -\frac{d}{dt} P_x(\tau_D(B) > t).$$

If we can show for fixed $x$ the series obtained from (A.2) by termwise differentiation with respect to $t$ converges uniformly for $t \geq T$, then by (A.5)

(A.8) $$f(t) = \sum_{j=1}^{\infty} B_j \left(\frac{|x|^2}{2t}\right)^{a_j/2} m_j\left(\frac{x}{|x|}\right) t^{-1}$$
$$\times \left[\frac{a_j}{2} {}_1F_1\left(\frac{a_j}{2}, a_j + \frac{n}{2}, -\frac{|x|^2}{2t}\right)\right.$$
$$\left. - \frac{a_j/2}{a_j + n/2} {}_1F_1\left(1 + \frac{a_j}{2}, 1 + a_j + \frac{n}{2}, -\frac{|x|^2}{2t}\right) \frac{|x|^2}{2t}\right].$$

Since ${}_1F_1(\alpha, \gamma, z) \to 1$ as $z \to 0$, the asymptotic (3.2) [with $p(D) = a_1/2$] is immediate. Thus, we need only prove uniform convergence in $t \geq T$ of the series in (A.8) for $x$ fixed. Choose $N$ so large that

(A.9) $$\frac{\Gamma(a/2 - n - 4)(2a + n/2)}{\Gamma(a/2)} \leq 2^{n+6}[a(a + n - 2)]^{-(n+3)/2} \quad \text{for } a \geq N.$$

Then choose $J$ so large that

$$a_j > N \vee 2n \quad \text{for } j \geq J.$$

This is possible by (A.1) and (A.3), since $\lambda_j \to \infty$ as $j \to \infty$ by (A.7).

By (A.4) and (A.6), for $t \geq T$ and $j \geq J$,

$$\left| B_j m_j\left(\frac{x}{|x|}\right) t^{-1} \right| \leq c(n) \lambda_j^{(n-1)/2} \frac{\Gamma((a_j + n)/2)}{\Gamma(a_j + n/2)} T^{-1}.$$

Hence, by Lemma A.2 with $k = n+4$, for $j \geq J$ and $t \geq T$, the $j$th term in the series (A.8) is bounded in absolute value by

$$T^{-1} c(n) \lambda_j^{(n-1)/2} \frac{\Gamma((a_j+n)/2)}{\Gamma(a_j+n/2)}$$

$$\times \left[ \frac{a_j}{2} \left( \frac{|x|^2}{2t} \right)^{n+4} \frac{\Gamma(a_j+n/2)\Gamma(a_j/2-n-4)}{\Gamma(a_j/2)\Gamma((a_j+n)/2)} \right.$$

$$\left. + \frac{1}{2} \left( \frac{|x|^2}{2t} \right)^{n+4} \frac{\Gamma(1+a_j+n/2)\Gamma(1+a_j/2-n-4)}{\Gamma(1+a_j/2)\Gamma((a_j+n)/2)} \right]$$

$$= \frac{T^{-1}}{2} c(n) \lambda_j^{(n-1)/2} \left( \frac{|x|^2}{2t} \right)^{n+4}$$

$$\times \left[ a_j \frac{\Gamma(a_j/2-n-4)}{\Gamma(a_j/2)} + \frac{(a_j+n/2)(a_j/2-n-4)\Gamma(a_j/2-n-4)}{(a_j/2)\Gamma(a_j/2)} \right]$$

$$\leq \frac{T^{-1}}{2} c(n) \lambda_j^{(n-1)/2} \left( \frac{|x|^2}{2t} \right)^{n+4} \frac{\Gamma(a_j/2-n-4)}{\Gamma(a_j/2)} \left[ 2a_j + \frac{n}{2} \right]$$

$$\leq \frac{T^{-1}}{2} c(n) \lambda_j^{(n-1)/2} \left( \frac{|x|^2}{2t} \right)^{n+4} 2^{n+6} [a_j(a_j+n-2)]^{-(n+3)/2} \qquad \text{by (A.9).}$$

Since $a_j(a_j+n-2) = \lambda_j$, this is the same as

$$\frac{T^{-1}}{2} c(n) 2^{n+6} \left( \frac{|x|^2}{2t} \right)^{n+4} \lambda_j^{-2} \leq c(n) T^{-1} 2^{n+5} \left( \frac{|x|^2}{2T} \right)^{n+4} \lambda_j^{-2}.$$

The uniform convergence for $t \geq T$ of the series in (A.8) follows from (A.7). □

LEMMA A.4. *For a bounded open set $D \subseteq \mathbb{R}^n$ with regular boundary,*

$$\log\left( \frac{d}{dt} P_x(\tau_D(B) \leq t) \right) \sim -\lambda_D t \qquad \text{as } t \to \infty,$$

*where $\lambda_D$ is the principal eigenvalue of $\frac{1}{2}\Delta$ on $D$ with Dirichlet boundary conditions.*

PROOF. There is an eigenfunction expansion

$$P_x(\tau_D(B) > t) = \sum_{k=1}^{\infty} e^{-\lambda_k t} \varphi_k(x) \int_D \varphi_k(y) \, dy,$$

where $0 < \lambda_1 \leq \lambda_2 \leq \cdots$ and $\{\varphi_k\}$ are the eigenvalues and eigenfunctions of $\frac{1}{2}\Delta$ on $D$ with Dirichlet boundary conditions. Moreover, for some $\alpha > 0$, for each $t_0 > 0$,

(A.10) $$|\varphi_k(x)|^2 \leq \alpha t_0^{-n/2} e^{\lambda_k t_0}$$

and

$$\sum_{k=1}^{\infty} e^{-\lambda_k t_0} < \infty$$

([45], pages 121–127). Then it is easy to show for $x \in D$,

$$\frac{d}{dt} P_x(\tau_D(B) \leq t) = -\frac{d}{dt} P_x(\tau_D(B) > t)$$
$$= \sum_{k=1}^{\infty} \lambda_k e^{-\lambda_k t} \varphi_k(x) \int_D \varphi_k(y) \, dy$$

and then that

$$\log\left(\frac{d}{dt} P_x(\tau_D(B) \leq t)\right) \sim -\lambda_1 t \qquad \text{as } t \to \infty. \qquad \square$$

**Acknowledgments.** I thank the referee for providing much simpler, elementary proofs of Lemmas 4.1 and A.1, as well as the names of the authors Allouba, Helms, Hochberg, Mądrecki, Nishioka, Krylov, Orsingher and Vanderbei cited in the Introduction.

DEPARTMENT OF MATHEMATICS
TEXAS A&M UNIVERSITY
COLLEGE STATION
TEXAS 77843-3368
USA
E-MAIL: deblass@math.tamu.edu